\documentclass[reqno,a4paper,12pt]{amsart}
\usepackage{amsmath,amsthm,amsfonts,amssymb, mathrsfs}
\usepackage{verbatim, graphicx, ifthen, enumitem}
\usepackage[T1]{fontenc}
\usepackage [applemac] {inputenc}
\usepackage{color, hyperref}

\allowdisplaybreaks

\let\oldmarginpar\marginpar
\renewcommand\marginpar[1]{\-\oldmarginpar[\raggedleft\footnotesize #1]%
{\raggedright\footnotesize #1}}
\newtheorem{theorem}{Theorem}[section]
\newtheorem*{theorem*}{Theorem}

\newtheorem{lemma}[theorem]{Lemma}

\newtheorem{proposition}[theorem]{Proposition}
\newtheorem{corollary}[theorem]{Corollary}
\theoremstyle{definition}

\newcommand{\Z}{\mathbb{Z}}

\newcommand{\abs}[1]{|#1|}

\newcommand{\inner}[2]{\left\langle #1,#2 \right\rangle}

\newcommand{\N}{\mathbb{d}}
\newcommand{\R}{\mathbb{R}}

\newcommand{\C}{\mathbb{C}}

\newcommand{\Bigoh}[1]{\mathcal{O} \left( #1 \right)}

\def\N{\mathbb{N}}
\def\Z{\mathbb{Z}}

\def\R{\mathbb{R}}
\def\R{\mathbb{R}}
\def\C{\mathbb{C}}

\def\1{\mathbf{1}}

\newcommand{\dif}{\mathrm{d}}
\newcommand{\e}{\mathrm{e}}
\newcommand{\im}{\mathrm{i}}
\newcommand{\norm}[1]{\|#1\|}

\renewcommand{\Re}{\operatorname{Re}}

\title{An operator theoretic approach to the prime number theorem}

\author[Olsen]{Jan-Fredrik Olsen} 
\address{Centre for Mathematical Sciences, Lund University, P.O. Box 118, SE-221 00 Lund, Sweden}
\email{jan-fredrik.olsen@math.lu.se}

\begin{document}

\maketitle

\begin{abstract}
	In this short note, we establish an operator theoretic version of the Wiener-Ikehara tauberian theorem, and point out how this leads to a new  proof of the Prime number theorem that should be accessible to anyone with  a basic knowledge of operator theory.
\end{abstract}

\section{Introduction}

We begin by stating a version of the Wiener-Ikehara tauberian theorem, due to Korevaar \cite{korevaar2005}. To this end, we remark that   the distributional Fourier transform of $v \in L^\infty(\R)$  such that $\lim_{\abs{x}\rightarrow \infty} v(x) = 0$  is called a \textit{pseudo-function}. Also, we denote the complex variable by $s = \sigma + \im t$.

\begin{theorem} 
	Let $S(t)$ be a non-decreasing function with support in $[0,\infty)$, and suppose that
	the Laplace transform
	\begin{equation*}
			 \mathcal{L}S(s) = 
		\int_0^\infty {S(u)} \e^{-su} \dif u
	\end{equation*}
	exists for $\sigma > 1$, and, for some constant $A$, let 
	\begin{equation*}
		g(s) = \mathcal{L}S(s) - \frac{A}{s-1}.
	\end{equation*}
	If $g(s)$ coincides with a pseudo-function on every bounded interval on the abscissa   $\sigma=1$ then
	\begin{equation*} \label{tauberian limit}
	 	\lim_{u\rightarrow \infty} \frac{S(u)}{\e^u} = A. 
	\end{equation*}
	Conversely, if this limit holds, then $g$ extends to a pseudo-function on $\sigma=1$.
\end{theorem}

We point out that Ikehara, a student of Wiener, originally established his tauberian theorem in order to find a simple analytic proof for the Prime number theorem. 

Below, we state and prove an   operator theoretic generalisation of the Wiener-Ikehara-Korevaar theorem. Since the machinery of operator theory allows us to avoid the delicate manipulations of limits required in Korevaar's proof, our hope is that this will provide a more accessible route to the Prime number theorem for anyone with a basic familiarity of operator theory. 

We also mention that J.-P.\,Kahane has a functional analytic  proof of the Prime number theorem, which is rather ingenious \cite{kahane96a}.

To motivate  our  approach,  we recall some ideas from \cite{olsen2010modzeta}. Specifically, we define, for intervals $I \subset \R$ symmetric with respect to the origin, the following operator on $L^2(I)$ (which we consider as a subspace of $L^2(\R)$):
\begin{align*}
	W_{I} : f  \longmapsto &\lim_{\epsilon \rightarrow 0^+} \frac{1}{\pi} \int_I f(\tau) \Re \frac{\zeta \big(1+\epsilon + \im(t-\tau)\big)}{1+\epsilon + \im(t-\tau)} \dif \tau,
\end{align*}
where $\zeta(s) = \sum_{n \in \N} n^{-s}$ is the  Riemann zeta function.    It is well-known that 
\begin{equation}  \label{central formula}
	\frac{\zeta(s)}{s} = \frac{1}{s-1} + \psi(s),
\end{equation}
where $\psi$ is an entire function.  Plugging  \eqref{central formula} into the formula for $W_I$, and noting that   the term $1/(s-1)$ leads to the appearance of the Poisson kernel, we   obtain   
\begin{equation} \label{aha}
	W_{I} f(t)
	= f(t)  + \Psi_{I}f,
\end{equation}
where $\Psi_I$ is readily seen to be a compact operator on $L^2(I)$.

Since     $\zeta(s)/s$ is the Laplace transform of $\pi_\N(\e^u)$, where   $\pi_\N$  denotes the counting function of the integers, it   follows from Plancherel's theorem that 
\begin{equation} \label{plancherelled version}
	W_{I}f(t) = \frac{1}{2\pi} \int_\R \frac{\pi_\N(\e^{\abs{u}})}{\e^{\abs{u}}}  \hat{f}(u) \e^{\im u t} \dif u.
\end{equation}
In particular, by the Fourier inversion formula, this means that
\begin{equation} \label{Psi-formula}
	\Psi_I f(t) = W_I f(t) - f(t) =\frac{1}{2\pi} \int_\R \Big( \frac{\pi_\N(\e^{\abs{u}})}{\e^{\abs{u}}} - 1 \Big) \hat{f}(u) \e^{\im u t} \dif u.
\end{equation}
This relation allows us to connect the operator theoretic properties of $W_I$ and $\Psi_I$ to the arithmetic properties of the counting function $\pi_\N$.  More generally, in \cite{olsen2010modzeta}, properties of the operators corresponding to   counting functions $\pi_K$ for  subsets $K \subset \N$ were studied.   In particular, it was observed that the resulting operator is closely related to  the frame operator for the set of vectors $\{ n^{\pm \im t-1/2}\}_{n \in K}$ which appears in the study   of Hilbert spaces of Dirichlet series (see also  \cite{olsen_saksman2012}).

We now state our main result, where   $\mathrm{Id}$ denotes the identity operator on $L^2(I)$.
\begin{theorem} \label{main theorem}
	Let $S(x)$ be a     non-decreasing function on $[0,\infty)$ and $I \subset \R$ an interval that is symmetric with respect to the origin.  Suppose that the Laplace transform $G(s) := \mathscr{L}\{S(\e^u)\}(s)$ exists for $\sigma > 1$, and for $\epsilon>0$ consider the operators
	\begin{align*}
		W_{S,I,\epsilon}: f \in L^2(I) \longmapsto   \frac{1}{\pi}  \int_I &f(\tau) \Re  {G \Big(1+\epsilon + \im(t-\tau)\Big)}  \dif \tau. 
	\end{align*}
	Then, for every $I$ sufficiently large, we have that as $\epsilon\rightarrow 0^+$, the operators $W_{S,I,\epsilon}$ converge in the weak operator norm   to  an operator $W_{S,I}$ satisfying
	$$W_{S,I} = A \mathrm{Id} + \Psi_{S,I},$$ for some constant $A$ and compact operator $\Psi_{S,I}$ on $L^2(I)$, if and only if $$\lim_{u \rightarrow \infty} \frac{S(\e^u)}{\e^u}  = A.$$
\end{theorem}

%

Before discussing the proof of the above result, we point out to the reader how the Prime number theorem follows from Theorem \ref{main theorem}. This argument should be clear to anyone familiar with Ikehara's theorem,  but we include it to keep the note self-contained.

\begin{corollary} \textbf{(The Prime number theorem)}
	Let $\pi_{\mathbb{P}}(x)$ be the counting function for the prime numbers. Then
	\begin{equation*}
		\lim_{x \rightarrow \infty} \pi_{\mathbb{P}}(x) \cdot \frac{\ln x}{x} = 1.
	\end{equation*}
\end{corollary}
\begin{proof}
	Let $\zeta_{\mathbb{P}}(s) = \sum_{p \text{ prime}} p^{-s}$. Then 
	\begin{equation} \label{prime formula 1}
		\frac{\zeta_\mathbb{P}(s)}{s} 
		 = \mathcal{L}\Big\{  {\pi_{\mathbb{P}}(\e^u)}  \Big\}(s).
	\end{equation}
	By taking the logarithm of the Euler product formula, and using a first order Taylor approximation on the terms $\log(1-p^{-s})$, we find that
	\begin{equation*}
		\log \zeta(s) = \zeta_{\mathbb{P}}(s) +  \sum_{\text{$p$ prime}} \Bigoh{p^{-2s}}.
	\end{equation*}
	 In combination  with \eqref{central formula}, this yields the well-known formula
	\begin{equation} \label{prime formula 2}
		\frac{\zeta_\mathbb{P}(s)}{s} = \log \frac{1}{s-1} + \psi_{\mathbb{P}}(s),
	\end{equation}
	where $\psi_{\mathbb{P}}(s)$ is analytic in a neighbourhood of $\{\Re s> 1\}$.
	
	We now combine formulas \eqref{prime formula 1} and \eqref{prime formula 2}, and differentiate, to obtain the relation
	\begin{equation} \label{ahahh}
		\mathcal{L}\big\{u\pi_{\mathbb{P}}(\e^{u}) \big\} = \frac{1}{s-1} - \psi_\mathbb{P}(\e^u).
	\end{equation}
	Finally,   by the same reasoning used to deduce \eqref{aha} from \eqref{central formula}, this implies that  for the choice $S(x) = \pi_{\mathbb{P}}(x)\ln x$, we have
	\begin{equation*}
		W_{S,I} f(t) = \mathrm{Id} + \Psi_{S,I}f,
	\end{equation*}
	where $\Psi_{S,I}$ is a compact operator on $L^2(I)$ for all bounded and symmetric intervals $I$. By Theorem \ref{main theorem},   the Prime number theorem now follows.
\end{proof}


\section{Proof of   theorem \ref{main theorem}}

We first suppose that ${S(e^u)}/{\e^u} \longrightarrow  A$ as $ {u}\rightarrow \infty$.
In the same way that we arrived at \eqref{Psi-formula}, we have
\begin{align*}
	\Psi_{S,I}f(t) 
	&= \frac{1}{2\pi} \int_\R \underbrace{\Big( \frac{S(\e^{\abs{u}})}{\e^{\abs{u}}} - A \Big)}_{:=h(u)} \hat{f}(u) \e^{\im u  t} \dif u.
\end{align*}
Since $\Psi_{S,I}$ is an operator on $L^2(I)$, where $I$ is a bounded interval, and $h(u)$ decays, it readily follows that $\Psi_{S,I}$ is compact (see, e.g., Lemma 2 in \cite{olsen2010modzeta}).

We now consider the converse implication (which is the one needed to prove the Prime number theorem).  To this end, let $S$ be a fixed non-decreasing function, and suppose that  $W_{S,I} = A \mathrm{Id} + \Psi_{S,I}$ for a compact operator $\Psi_{S,I}$ for $I$ large enough.  

First, we show that $\Psi_{S,I}$ bounded implies that $g(u) : = S(\e^u)/\e^u$ is bounded. 
Letting $\{e_n\}_{n \in \Z}$ denote the standard orthonormal exponential basis for $L^2(I)$, we compute
\begin{equation*} 
	\begin{aligned}
	\inner{W_{S,I,\epsilon} e_n}{e_n } 
	=  \frac{1}{2\pi} \int_\R g(u)  \abs{I}  \Big(\frac{\sin( u\abs{I}/2 )  }{ u\abs{I}/2 - \pi n}\Big)^2  \e^{-\epsilon \abs{u}} \dif u.
	\end{aligned}
\end{equation*}
To obtain a contradiction, suppose   that $g(u)$ is not bounded. Then there exists a sequence $u_k$ of positive numbers tending to infinity, such that $g(u_k)\geq k$. But since $S$ is non-decreasing, we get, for $\Delta u>0$, that
\begin{equation*}
	g(u_k + \Delta u) = \frac{S(\e^{u_k + \Delta u})}{\e^{u_k + \Delta u}} \geq \frac{S(\e^{u_k})}{\e^{u_k }} \frac{1}{\e^{\Delta u}} \geq \frac{k}{\e^{\Delta u}}.
\end{equation*}
It now follows   that, for sufficiently large $I$,   the sequence $\inner{W_{S,I}e_n}{e_n}$ is unbounded, which is absurd.

%

Next, we recall that if $\Psi_{S,I}$ is compact then 
\begin{equation*}
	\inner{\Psi_{S,I}e_n}{e_n} \rightarrow 0 \quad \text{as } n \rightarrow \infty. 
\end{equation*}
Using the fact that $g(u)$, and therefore $h(u)$, defined as in the first part of the proof, is bounded,   a straight-forward computation gives
\begin{equation} \label{computation}
	\begin{aligned}
	\inner{\Psi_{S,I} e_n}{e_n } 
	&= \frac{1}{2\pi}  \int_\R h(u)  \abs{I}  \Big(\frac{\sin( u\abs{I}/2 )  }{ u\abs{I}/2 - \pi n}\Big)^2  \dif u.
	\end{aligned}
\end{equation}
To arrive at a contradiction, suppose that  the limit  $h(u) \rightarrow 0$
does not hold as $\abs{u}\rightarrow \infty$. There are now two cases. In the first, we suppose  that there exists an $\epsilon>0$ and an unbounded sequence $u_k$ of real numbers so that
\begin{equation*}
	h(u_k) = \frac{S(\e^{ \abs{u_k}})}{\e^{ \abs{u_k}}} - A  \geq \epsilon.
\end{equation*}
As above, it follows that 
there exists a fixed $\Delta u >0$ so that for all $k \in \N$ and  $u \in [u_k, u_k + \Delta u]$ we have
\begin{equation*}
	h(u)  \geq \frac{\epsilon}{2}.
\end{equation*}
It is now straight-forward to apply this estimate to the integral expression in \eqref{computation} to see that for all $I$ large enough, there exists a constant $c = c(I)>0$ so that for  infinitely many $n$, we have
\begin{equation*}
	\abs{\inner{\Psi_{S,I}e_n}{e_n}} \geq c.
\end{equation*}

In the remaining case, we suppose that there exist an unbounded sequence of  real numbers $u_k$ so that $h(u) \leq -\epsilon$.
This case is settled as above with the adjustment that we consider intervals of the type   $[u_k - \Delta u, u_k]$. \qed

\section{Remarks}

While working on a manuscript containing these results, the author was made aware of similar efforts by Franz Luef and Eirik Skrettingland, who, in a rather technical paper, obtain a generalised version of Wiener's tauberian theorem for operators (see \cite{luef_skrettingland2020}). Curiously enough, the results of Luef and Skrettingland only apply when the   underlying space is $L^2(\R^n)$ for even $n$, and therefore do not seem to imply the results in this note. We plan to investigate the relation between these Tauberian theorems in a forthcoming publication. 


%

\bibliographystyle{amsplain}

\def\cprime{$'$} \def\cprime{$'$} \def\cprime{$'$} \def\cprime{$'$}
\providecommand{\bysame}{\leavevmode\hbox to3em{\hrulefill}\thinspace}
\providecommand{\MR}{\relax\ifhmode\unskip\space\fi MR }
\providecommand{\MRhref}[2]{%
  \href{http://www.ams.org/mathscinet-getitem?mr=#1}{#2}
}
\providecommand{\href}[2]{#2}

\end{document}